\documentclass[12pt]{article}

\usepackage{a4wide}

\newcommand{\Ga}{\Gamma}

     \title{A proof of Kummer's theorem}

     \author{Bruno Gauthier
             \thanks{ \texttt{Bruno.Gauthier@univ-mlv.fr},
             \textit{Institut Gaspard Monge}, 
             Universit\'e Marne-la-Vall\'ee, Cit\'e Descartes, 5 bd Descartes,
             Champs sur marne, F-77454 Marne-la-Vall\'ee Cedex 2.}
            }

\begin{document}

     \maketitle

     \begin{abstract}
     Following suggestions of T.~H.~Koornwinder \cite{b-tom}, we
     give a new proof of Kummer's theorem involving Zeilberger's
     algorithm, the WZ method and asymptotic estimates.
     In the first section, we recall a classical proof given by
     L.~J.~Slater \cite{b-slater}. The second section discusses
     the new proof, in the third section sketches of similar proofs for
     Bailey's and Dixon's theorems are given.
     \end{abstract}

%\begin{keyword}
%Nonterminating hypergeometric series; Kummer's theorem, Bailey's theorem,
%Dixon's theorem, Zeilberger's algorithm.
%\end{keyword}

The author is grateful to Peter Paule for his helpful comments.

\section{Slater's proof}

In classical hypergeometric theory, formulae concerning summation,
transformation and contiguous relations are provided in order to
manipulate expressions that represent special functions.
One of those is Kummer's non-terminating summation theorem:
\begin{equation} \label{e:kummerth}
 {}_2F_1\left[{a,b \atop 1+a-b};-1\right]
     = \frac{ \Ga(1+\frac{a}{2}) \: \Ga(1+a-b) }
            { \Ga(1+a) \: \Ga(1+\frac{a}{2}-b) } \ \mbox{,}
\end{equation}
% $ Re(1+a-b-a-b) > -1 $,
in which $ Re(b) < 1 $ ensures the series to be convergent.
%in order that the series is convergent.

The first step in Slater's proof is to prove the following quadratic
transformation, due to Kummer:
\begin{equation} \label{e:kummerqua}
 {}_2F_1\left[{a,b \atop 1+a-b};z\right]
     = (1-z)^{-a}
       \: {}_2F_1\left[{\frac{a}{2},\frac{1}{2}+\frac{a}{2}-b
                       \atop 1+a-b};\frac{-4\:z}{(1-z)^2}\right] \ \mbox{.}
\end{equation}
The proof of this transformation is as follows.\\
The right-hand side of (\ref{e:kummerqua}) can be written:
\begin{eqnarray*} \displaystyle
& & (1-z)^{-a} \: \sum_{k=0}^\infty 
               \frac{(\frac{a}{2})_k \: (\frac{1}{2}+\frac{a}{2}-b)_k}
                    {(1+a-b)_k \: k!}
               \: \frac{(-4)^k \: z^k}{(1-z)^{2\:k}}\\
&=& \sum_{k=0}^\infty 
    \frac{(\frac{a}{2})_k \: (\frac{1}{2}+\frac{a}{2}-b)_k \: (-4)^k}
         {(1+a-b)_k \: k!}
    \: z^k \: (1-z)^{-a-2\:k} \ \mbox{,}\\
&=& \sum_{k=0}^\infty 
    \frac{(\frac{a}{2})_k \: (\frac{1}{2}+\frac{a}{2}-b)_k \: (-4)^k}
         {(1+a-b)_k \: k!}
    \: z^k \: {}_1F_0\left[{a+2\:k \atop -};z\right] 
              %\sum_{n=0}^\infty \frac{(a+2\:k)_n}{n!} \: z^n 
                                   \quad \mbox{by the binomial theorem,}\\
&=& \sum_{k=0}^\infty \sum_{n=0}^\infty 
    \frac{(\frac{a}{2})_k \: (\frac{1}{2}+\frac{a}{2}-b)_k \: (-4)^k}
         {(1+a-b)_k \: k!}
    \: \frac{(a+2\:k)_n}{n!} \: z^{k+n} \ \mbox{.}
\end{eqnarray*}
Extracting the coefficient of $z^N$ yields:
\begin{eqnarray}
 %\displaystyle
 &   & \sum_{k=0}^{N}
    \frac{(\frac{a}{2})_k \: (\frac{1}{2}+\frac{a}{2}-b)_k
          \: (-4)^k \: (a+2\:k)_{N-k}} 
         {(1+a-b)_k \: (N-k)! \: k!} \nonumber \\
\qquad & = & \frac{(a)_N}{N!}
    \: \sum_{k=0}^{N}
       \frac{(\frac{1}{2}+\frac{a}{2}-b)_k \: (a+N)_k \: (-N)_k}
            {(1+a-b)_k \: (\frac{1}{2}+\frac{a}{2})_k \: k!}
 \label{e:beforesaal}
\end{eqnarray}
since:
\[ \frac{1}{(N-k)!} = \frac{(-1)^k \: (-N)_k}{N!} \]
and
\[ (a+2\:k)_{N-k} = \frac{(a)_{N+k}}{(a)_{2\:k}}
                  = \frac{(a)_N \: (a+N)_k}{4^k \: (\frac{a}{2})_k
                    \: (\frac{1}{2}+\frac{a}{2})_k} \ \mbox{.} \]
The right-hand side of (\ref{e:beforesaal}) is summable by
Saalsch{\" u}tz theorem \cite[(III.2)]{b-slater}:
\[
   \frac{(a)_N}{N!}
   \: \frac{(\frac{1}{2}+\frac{a}{2})_N \: (1-b-N)_N}
           {(1+a-b)_N \: (\frac{1}{2}-\frac{a}{2}-N)_N} \ \mbox{.}
\]
As this is the coefficient of $z^N$ in the right-hand side
of (\ref{e:kummerqua}), we obtain:
\begin{eqnarray*}
   \sum_{N=0}^{\infty} \frac{(a)_N}{N!}
      \: \frac{(\frac{1}{2}+\frac{a}{2})_N \: (1-b-N)_N}
              {(1+a-b)_N \: (\frac{1}{2}-\frac{a}{2}-N)_N} \: z^N
  & = & \sum_{N=0}^{\infty} \frac{(a)_N \: (b)_N}{(1+a-b)_N \: N!} \: z^N \\
  & = & {}_2F_1\left[{a,b \atop 1+a-b};z\right] \ \mbox{,}
\end{eqnarray*}
which proves Kummer's quadratic transformation.\\

The second step is to set $z$ to $-1$ in (\ref{e:kummerqua}). Therefore
\begin{equation} \label{e:kummerquaspe}
 {}_2F_1\left[{a,b \atop 1+a-b};-1\right]
     = 2^{-a}
       \: {}_2F_1\left[{\frac{a}{2},\frac{1}{2}+\frac{a}{2}-b
                       \atop 1+a-b};1\right] \ \mbox{.}
\end{equation}
Note that these series are convergent if $ Re(b) < 1 $.\\
We can sum the series on the right-hand side
by Gauss's theorem \cite[(III.3)]{b-slater}. This yields
\begin{equation}
 {}_2F_1\left[{a,b \atop 1+a-b};-1\right]
 = 2^{-a}
   \: \frac{\Ga(1+a-b) \: \Ga(\frac{1}{2})}
           {\Ga(1+\frac{a}{2}-b) \: \Ga(\frac{1}{2}+\frac{a}{2})} \ \mbox{.}
\end{equation}
Since $ \Ga(\frac{1}{2}+\frac{a}{2}) \: \Ga(\frac{a}{2}+1) 
             = 2^{a} \: \Ga(\frac{1}{2}) \: \Ga(a+1) $,
we finally obtain Kummer's theorem:
\begin{equation}
 {}_2F_1\left[{a,b \atop 1+a-b};-1\right]
     = \frac{ \Ga(1+a-b) \: \Ga(1+\frac{a}{2})}
            { \Ga(1+\frac{a}{2}-b) \: \Ga(1+a)} \ \mbox{.}
\end{equation}

%\newpage

\section{A new proof}

In formula (\ref{e:kummerth}), let us replace $a$ by $a+2\:n$
where $n$ is a new variable that stands for a nonnegative integer.
The resulting formula
\begin{equation} \label{e:kummervar}
 {}_2F_1\left[{a+2\:n,b \atop 1+a+2\:n-b};-1\right]
     = \frac{ \Ga(1+\frac{a}{2}+n) \: \Ga(1+a+2\:n-b) }
            { \Ga(1+a+2\:n) \: \Ga(1+\frac{a}{2}+n-b) }
\end{equation}
is clearly equivalent to (\ref{e:kummerth}).\\
The introduction of the free parameter $n$ is the key to this new
computer-assisted proof.

Formula (\ref{e:kummervar}) can be written as:
\begin{equation} \label{e:p1}
    \frac{\sum_k f(n,k)}{S(n)} = 1 \ \mbox{,}
\end{equation}
where \[ \displaystyle f(n,k) = \frac{(a+2\:n)_k \: (b)_k}{(1+a+2\:n-b)_k}
                  \: \frac{(-1)^k}{k!} \]
and   \[ \displaystyle S(n) = \frac{\Ga(1+\frac{a}{2}+n) \: \Ga(1+a+2\:n-b)}
                    {\Ga(1+a+2\:n) \: \Ga(1+\frac{a}{2}+n-b)} \ \mbox{.} \]
(The coefficient $2$ in the substitution $a \rightarrow a+2\:n$ makes $S(n)$
hypergeometric in $n$.)\\
Formula (\ref{e:p1}) can be rewritten as:
\begin{equation} \label{e:p2}
\sum_k F(n,k) = 1 \qquad \mbox{where }
F(n,k) = \frac{f(n,k)}{S(n)} \mbox{.}
\end{equation}
We shall actually prove this last formula. Using Gauthier's Maple
package \texttt{HYPERG} \cite{b-user}, we apply Zeilberger's algorithm
to obtain:
\begin{equation} \label{e:kummerWZ}
F(n,k) - F(n+1,k) = G(n,k+1) - G(n,k) \ \mbox{,}
\end{equation}
with $ G(n,k) = F(n,k) \: C(n,k) $ and $ \displaystyle C(n,k)
= -\frac{(b-1)\:k}{(1+a+2\:n-b+k) \: (a+2\:n)} $.\\
This $C(n,k)$ is the so-called certificate of the WZ-pair $(F,G)$.\\
After summing (\ref{e:kummerWZ}) for $ k $ from $ 0 $ to $K$,
the right-hand side telescopes:
\begin{eqnarray}
 \sum_{k=0}^{K} F(n,k) - \sum_{k=0}^{K} F(n+1,k)  
 & = & \sum_{k=0}^{K} G(n,k+1) - \sum_{k=0}^{K} G(n,k) \\
 & = & G(n,K+1) - G(n,0) \\
 \label{e:WZsum}
 & = & F(n,K+1) \: C(n,K+1) \ \mbox{.}
\end{eqnarray} 
We now let $K$ tend to infinity. Clearly,
\[ \lim_{K \rightarrow +\infty} C(n,K+1) = -\frac{b-1}{a+2\:n} \ \mbox{.} \]
To determine $ \displaystyle \lim_{K\rightarrow+\infty} F(n,K+1) $, we use
the well-known estimate (see e.g. \cite{b-tom}):
       \begin{equation} \label{e:lim}
       \frac{\Ga(a+k)}{\Ga(b+k)} \sim k^{a-b}
       \quad \mbox{as } k \rightarrow \infty \ \mbox{.}
       \end{equation}
First,
\begin{eqnarray*}
 F(n,k) & = & \frac{(a+2\:n)_k \: (b)_k}{(1+a+2\:n-b)_k}
              \: \frac{(-1)^k}{k!} \: \frac{1}{S(n)} \ \mbox{,} \\
        & = & \frac{\Ga(a+2\:n+k) \: \Ga(b+k) \: (-1)^k}
                   {\Ga(1+a+2\:n-b+k) \: \Ga(k+1)}
              \: T(n) \ \mbox{,}
\end{eqnarray*}
where
$\displaystyle T(n)=\frac{(a+2\:n)\:\Ga(1+\frac{a}{2}+n-b)}
                         {\Ga(b)\:\Ga(1+\frac{a}{2}+n)}$
does not depend on $k$.\\
Then,
\begin{eqnarray*}
 |F(n,k)| & \sim & k^{a+2n-(1+a+2n-b)} \: k^{b-1} \: T(n) \ \mbox{,} \\
          & \sim & k^{2b-2} \: T(n) \ \mbox{.}
\end{eqnarray*}
We know that $ Re(b)<1 $ (necessary condition for the convergence
of the series in (\ref{e:kummerth})).\\
Therefore
\[ \lim_{K \rightarrow +\infty} F(n,K+1) = 0 \ \mbox{.} \]

Expressing as before $F(n,k)$ in terms of Gamma functions and using the same
aymptotics, it is easy to find the limit of $F(n,k)$ when $n$ tends to
infinity:
\begin{eqnarray*}
  F(n,k)  & = & \frac{(a+2\:n)_k \: (b)_k}{(1+a+2\:n-b)_k}
                \: \frac{(-1)^k}{k!}
                \: \frac{\Ga(1+a+2\:n) \: \Ga(1+\frac{a}{2}+n-b)}
                       {\Ga(1+\frac{a}{2}+n) \: \Ga(1+a+2\:n-b)} \\
          & = & \frac{\Ga(a+2\:n+k) \: \Ga(1+a+2\:n-b)}
                     {\Ga(a+2\:n)\:\Ga(1+a+2\:n-b+k)}
                \: \frac{(b)_k \: (-1)^k}{k!} \\
          &   & \times \frac{\Ga(1+a+2\:n) \: \Ga(1+\frac{a}{2}+n-b)}
                        {\Ga(1+a+2\:n-b) \: \Ga(1+\frac{a}{2}+n)} \\
 |F(n,k)| &\sim& \left| (2\:n)^{k} \: (2\:n)^{-k} \: \frac{(b)_k \: (-1)^k}{k!}
                \: (2\:n)^b \: n^{-b} \right| \ \mbox{.}
\end{eqnarray*}
So,
\begin{equation}
\lim_{n \rightarrow +\infty} F(n,k)
   = \frac{(b)_k \: (-1)^k}{k!} \: 2^b \ \mbox{.}
\end{equation}

By (\ref{e:lim}), we have $ 1/S(n) \sim 2^b $ when $n$ tends to infinity, so
for all $n \geq n_0$, $1/S(n) < A$ where $A$ is a positive constant.
Therefore,
\begin{equation} \label{e:convS}
 |F(n,k)| = \left|\frac{f(n,k)}{S(n)}\right| \leq A \: |f(n,k)| \ \mbox{.}
\end{equation}
Moreover,
\begin{equation} \label{e:majorF1}
 |f(n,k)| =    \frac{ |(a+2\:n)_k| \: |(b)_k| }{ |(1+a+2\:n-b)_k| \: k! }
          \leq \frac{ |(b)_k| }{ k! } \ \mbox{,}
\end{equation}
because
\[
  \forall i \geq 0 \ \mbox{,} \qquad
  \frac{|a+2\:n+i|}{|1+a+2\:n-b+i|} < 1 
  \quad \mbox{when } n \mbox{ is big enough,}
\]
and
\begin{eqnarray*}
  &      & \frac{ |(a+2\:n)_k| }{ |(1+a+2\:n-b)_k| } \\[0.2cm]
  &   =  & \frac{ | (a+2\:n) \: (a+2\:n+1) \ldots (a+2\:n+k-1) |}
                { | (1+a+2\:n-b) \: (1+a+2\:n-b+1) \ldots (1+a+2\:n+b+k-1) |}
           \ \mbox{,} \\[0.2cm]
  & \leq & \frac{ |a+2\:n| \: |a+2\:n+1| \ldots |a+2\:n+k-1| }
                { |1+a+2\:n-b| \: |1+a+2\:n-b+1| \ldots |1+a+2\:n+b+k-1| }
           \ \mbox{,} \\[0.2cm]
  & \leq & 1 \ \mbox{.}
\end{eqnarray*}
We can finally establish (use (\ref{e:convS}), (\ref{e:majorF1})
and (\ref{e:lim})) that
\begin{equation} \label{e:majorF2}
 |F(n,k)| \leq A \: \left| \frac{(b)_k}{k!} \right|
                  \sim A \: \left| \frac{k^{b-1}}{\Ga(b)} \right|  \mbox{.}
\end{equation}
%but this point doesn't prove the convergence of the series yet.

%  We define a new sequence $ H(n,k) := F(n,2\:k) + F(n,2\:k+1) $.
%  We have:
%  \begin{eqnarray*}
%  |H(n,k)| & = & |F(n,2\:k)|
%               \: \left| \frac{(4-4\:b)\:k+n\:(2-2\:b)+1+a-b-a\:n}
%                              {4\:k^2+(4+4\:n-2\:b+2\:a)\:k+2\:n+a+1-b} \right|
%                                                \ \mbox{,} \\
%           & \leq & |F(n,2\:k)| \: \frac{1}{k} \: B \ \mbox{,} \\
%           & \leq & C \: \frac{|(b)_{2\:k}|}{(2\:k)!} \: \frac{1}{k} \ \mbox{,}
%  %         & \leq & D \: \frac{|(b)_{2\:k}|}{(k+1)!} \ \mbox{,}
%  \end{eqnarray*}
%  This last expression is equivalent to $ D \: k^{Re(b)-2} $
%  (where $B$, $C$ and $D$ are positive constants).
%  It is the general term of a convergent series, so $\sum_{k \geq 0} H(n,k)$
%  is dominated by a convergent series.
%  %Moreover $ \lim_{k\rightarrow +\infty} |F(n,k)| = 0 $, so
%  %we can establish that $\sum_{k \geq 0} F(n,k)$ converges.

For $Re(b)<0$, this last expression is the general term of a convergent series,
so $\sum_{k \geq 0} F(n,k) $ is dominated by a convergent series.

The limit of formula (\ref{e:WZsum}) when $K$ tends to infinity gives
\[
   \sum_{k=0}^{\infty} F(n,k) - \sum_{k=0}^{\infty} F(n+1,k) = 0 \ \mbox{.}
\]
Hence, the sum $\sum_{k=0}^{\infty} F(n,k) $ is independent of $n$. Then
\begin{eqnarray*}
\sum_{k=0}^{\infty} F(n,k)
   & = & \lim_{n \rightarrow +\infty} \sum_{k=0}^{\infty} F(n,k) \ \mbox{,} \\
%   & = & \lim_{n \rightarrow +\infty}
%              \sum_{k=0}^{\infty} (F(n,2\:k)+F(n,2\:k+1)) \ \mbox{,} \\
%   & = & \lim_{n \rightarrow +\infty} \sum_{k=0}^{\infty} H(n,k) 
   & = & \sum_{k=0}^{\infty} \lim_{n \rightarrow +\infty} F(n,k) \ \mbox{,}
\end{eqnarray*}
where the last equality is justified by dominated convergence.
Moreover,
%\begin{eqnarray*}
%\sum_{k=0}^{\infty} F(n,k)
%   & = & \sum_{k=0}^{\infty} 
%         \lim_{n \rightarrow +\infty} (F(n,2\:k)+F(n,2\:k+1)) \ \mbox{,} \\
%   & = & \sum_{k=0}^{\infty} 2^b \:
%           \left(  \frac{(b)_{2\:k} \: (-1)^{2\:k}}{(2\:k)!} 
%            + \frac{(b)_{2\:k+1} \: (-1)^{2\:k+1}}{(2\:k+1)!} \right) \\
%   & = & 2^b \: \sum_{k=0}^{\infty} \: \frac{(b)_{k} \: (-1)^k}{k!}
%         \ \mbox{,} \\
%   & = & 2^b \: (1-(-1))^b = 1 \ \mbox{,} \qquad \mbox{by binomial theorem.}
%\end{eqnarray*}

\begin{eqnarray*}
\sum_{k=0}^{\infty} F(n,k)
   & = & \sum_{k=0}^{\infty} 
         \lim_{n \rightarrow +\infty} (F(n,k) \ \mbox{,} \\
   & = & \sum_{k=0}^{\infty} 2^b \:  \frac{(b)_k \: (-1)^k}{(k)!} \\
   & = & 2^b \: \sum_{k=0}^{\infty} \: \frac{(b)_{k} \: (-1)^k}{k!}
         \ \mbox{,} \\
   & = & 2^b \: (1-(-1))^b = 1 \ \mbox{,} \qquad \mbox{by binomial theorem.}
\end{eqnarray*}
This is precisely formula (\ref{e:p2}). So, we proved ---~with the restriction $
Re(b)<0 $~--- the validity of Kummer's theorem
(specializing $n$ by $0$ in (\ref{e:kummervar})).

Kummer's theorem (\ref{e:kummerth}) can be rewritten as:
\[
\sum_k M(a,b,k) = \frac{ \Ga(1+\frac{a}{2}) \: \Ga(1+a-b) }
                       { \Ga(1+a) \: \Ga(1+\frac{a}{2}-b) } \ \mbox{,}
\]
with
$\displaystyle M(a,b,k) = \frac{(a)_k \: (b)_k}{(1+a-b)_k \: k!} \: (-1)^k $.\\
We derive by Zeilberger's algorithm, a recurrence for $M$ with respect
to the parameter $b$:
\[
    (a-2\:b) \: M(a,b,k) + (-2\:a+2\:b) \: M(a,b+1,k)
  = G'(a,b,k+1) - G'(a,b,k)
\]
with $ G'(a,b,k) = M(a,b,k) \: \frac{(a-b+k)\:k}{b} $ and $ Re(b)<0 $.
After summing both sides with respect to $k$ from $0$ to $K-1$, it gives
\[
    (a-2\:b) \: \sum_{k=0}^{K-1} M(a,b,k)
  + (-2\:a+2\:b) \: \sum_{k=0}^{K-1} M(a,b+1,k)
  = G'(a,b,K) - G'(a,b,0) \ \mbox{,}
\]
and then we let $K$ tend to infinity to obtain
\[
    (a-2\:b) \: \sum_{k=0}^{\infty} M(a,b,k)
  + (-2\:a+2\:b) \: \sum_{k=0}^{\infty} M(a,b+1,k)
  = 0 \ \mbox{,}
\]  
because $ G'(a,b,0) = 0 $ and
$ \displaystyle \lim_{K \rightarrow +\infty} G'(a,b,K) = 0 $.

We have now for any $b$ with $Re(b)<0$
\begin{eqnarray*}
  \sum_{k=0}^{\infty} M(a,b+1,k)
  & = & \frac{\frac{a}{2}-b}{a-b} \: \sum_{k=0}^{\infty} M(a,b,k) \ \mbox{,} \\
  & = & \frac{\frac{a}{2}-b}{a-b}
        \: \frac{ \Ga(1+\frac{a}{2}) \: \Ga(1+a-b) }
                { \Ga(1+a) \: \Ga(1+\frac{a}{2}-b) } \ \mbox{,} \\
  & = & \frac{ \Ga(1+\frac{a}{2}) \: \Ga(a-b) }
             { \Ga(1+a) \: \Ga(\frac{a}{2}-b) } \ \mbox{.}
\end{eqnarray*}
Let us replace $b+1$ by $B$, we finally conclude that
\[
   {}_2F_1\left[{a,B \atop 1+a-B};-1\right]
   = \frac{ \Ga(1+\frac{a}{2}) \: \Ga(1-a-B) }
          { \Ga(1+a) \: \Ga(1+\frac{a}{2}-B) } \ \mbox{,}
\]
for all values of $B$ such that $Re(B)<1$.

\textbf{Remark:} We can also derive this result directly
from (\ref{e:majorF2}) by a convergence-acceleration argument.
For this, we define a new sequence $ H(n,k) := F(n,2\:k) + F(n,2\:k+1) $,
and we find that $\sum_{k \geq 0} H(n,k)$ is dominated by a convergent series.
Hence, we have
\begin{eqnarray*}
\sum_{k=0}^{\infty} F(n,k)
   & = & \lim_{n \rightarrow +\infty} \sum_{k=0}^{\infty} F(n,k) \ \mbox{,} \\
   & = & \lim_{n \rightarrow +\infty} 
         \sum_{k=0}^{\infty} (F(n,2\:k)+F(n,2\:k+1)) \ \mbox{,} \\
   & = & \lim_{n \rightarrow +\infty} \sum_{k=0}^{\infty} H(n,k) 
     =   \sum_{k=0}^{\infty} \lim_{n \rightarrow +\infty}  H(n,k) \ \mbox{,} \\
   & = & \sum_{k=0}^{\infty}
         \lim_{n \rightarrow +\infty} (F(n,2\:k)+F(n,2\:k+1)) \ \mbox{,} \\
   & = & \sum_{k=0}^{\infty} 2^b \:
           \left(  \frac{(b)_{2\:k} \: (-1)^{2\:k}}{(2\:k)!}
            + \frac{(b)_{2\:k+1} \: (-1)^{2\:k+1}}{(2\:k+1)!} \right) \\
   & = & 2^b \: \sum_{k=0}^{\infty} \: \frac{(b)_{k} \: (-1)^k}{k!}
         \ \mbox{,} \\
   & = & 1 \ \mbox{.}
\end{eqnarray*}

\section{Other summation theorems}

The combination of Zeilberger's algorithm and asymptotic estimates
suggests an approach to automatically prove identities involving
hypergeometric series. It is well suited to summation theorems.
As an example, we now give two apparently new proofs of classical
hypergeometric identities then a list of summation theorems
that the package \texttt{HYPERG} can prove automatically.

\subsection{Bailey's theorem}

\begin{equation} \label{e:baileyth}
{}_2F_1\left[{a,1-a \atop b};\frac{1}{2}\right]
       = \frac{ \Ga(\frac{b}{2}) \: \Ga(\frac{1+b}{2}) }
              { \Ga(\frac{a+b}{2}) \: \Ga(\frac{1-a+b}{2}) } \ \mbox{,}
\end{equation}
Following the same scheme, let $ b \rightarrow b+2\:n $.
Formula (\ref{e:baileyth}) can be written:
\begin{equation} \label{e:bai1}
\sum_k F(n,k) = 1 \ \mbox{,}
\end{equation}
where
$ \displaystyle F(n,k) = \frac{f(n,k)}{S(n)} $, \\
with 
$ \displaystyle
f(n,k) = \sum_k \frac{(a)_k \: (1-a)_k}{(b+2\:n)_k}
        \: \left(\frac{1}{2}\right)^k \quad \mbox{and} \quad 
S(n) = \frac{\Ga(\frac{b}{2}+n) \: \Ga(\frac{1}{2}+\frac{b}{2}+n)}
            {\Ga(\frac{a+b}{2}+n) \: \Ga(\frac{1-a+b}{2}+n)} \ \mbox{.} $\\
Zeilberger's algorithm applied to $F(n,k)$ succeeds to finding a WZ-pair:
\begin{equation} \label{e:baileyWZ}
F(n,k) - F(n+1,k) = G(n,k+1) - G(n,k) \ \mbox{,}
\end{equation}
where
$ \displaystyle G(n,k) = F(n,k) \: \left(-2\:\frac{k}{b+2\:n+k}\right) $.\\
As with Kummer's theorem, we have:
\begin{eqnarray}
G(n,0)                              & = & 0 \ \mbox{,} \\
\lim_{k \rightarrow +\infty} G(n,k) & = & 0 \ \mbox{,}
\end{eqnarray} 
and so, $ \displaystyle \sum_k F(n,k) $ is independent of $n$.\\
In this case, $ \displaystyle \lim_{n \rightarrow +\infty} F(n,k)
= \delta_{k,0} $, so that:
\begin{equation}
\sum_{k=0}^{\infty} F(n,k)
  = \lim_{n \rightarrow +\infty} \sum_{k=0}^{\infty} F(n,k) 
  = \sum_{k=0}^{\infty} \lim_{n \rightarrow +\infty} F(n,k)
  = \sum_{k=0}^{\infty} \delta_{k,0} = 1 \ \mbox{,}
\end{equation}
which proves Bailey's theorem.

\subsection{Dixon's theorem}

\begin{eqnarray} \label{e:dixonth}
&   & {}_3F_2\left[{a,b,c \atop 1+a-b,1+a-c};1\right] \nonumber \\
& = & \frac{ \Ga(1+a-b) \: \Ga(1+a-c) \: \Ga(1+\frac{a}{2})
                \: \Ga(1+\frac{a}{2}-b-c) }
              { \Ga(1+a) \: \Ga(1+\frac{a}{2}-b)
                \: \Ga(1+\frac{a}{2}-c) \: \Ga(1+a-b-c) } \ \mbox{,}
\end{eqnarray}
where $ Re(2+a-2\:b-2\:c) > 0 $.\\
Once again, we replace $ a $ by $a+2\:n$. With the same notation
(\ref{e:bai1}), we now have:
\begin{eqnarray*}
f(n,k) & = & \sum_k \frac{(a+2n)_k \: (b)_k \: (c)_k}
                     {(1+a+2n-b)_k \: (1+a+2n-c)_k \: k!} \ \mbox{,} \\
        \mbox{and} \quad
S(n) & = & \frac{\Ga(1+a+2n-b) \Ga(1+a+2n-c)
             \Ga(1+\frac{a}{2}+n) \Ga(1+\frac{a}{2}+n-b-c)}
            {\Ga(1+a+2n) \Ga(1+\frac{a}{2}+n-b)
             \Ga(1+\frac{a}{2}+n-c) \Ga(1+a+2n-b-c)} \ \mbox{.}
\end{eqnarray*}
The WZ-pair given by Zeilberger's algorithm is
\begin{equation} \label{e:dixonWZ}
F(n,k) - F(n+1,k) = G(n,k+1) - G(n,k) \ \mbox{,}
\end{equation}
with
$ \displaystyle G(n,k) = F(n,k) \: C(n,k) $, where the computed
certificat $C(n,k)$ is:
{ \footnotesize
\begin{verbatim}
-(-2-8*a*n*c-a+4*b+4*c-2*n-2*n*b*c-2*a*b*k-2*a*k*c-a*b*c-8*a*n*b
-4*n*k*c-4*n*b*k+a*b+2*n*b-2*a^2*c-8*n^2*c+a*c+2*n*c-2*b^2+3*a^2
+12*a*n+12*n^2+12*a^2*n-2*a^2*b+24*a*n^2+a*k^2-8*n^2*b+2*n*k^2
+2*a^3+16*n^3-6*c*b+3*k*a^2+3*k*a+6*k*n+12*k*n^2+12*k*a*n+2*b^2*c
+2*b*c^2-2*c^2)*k/(2+a+2*n-2*b-2*c)/(1+a+2*n-c+k)/(1+a+2*n-b+k)
/(a+2*n)
\end{verbatim} }
\noindent It is easy to establish the following limits:
\begin{eqnarray*}
G(n,0)                              & = & 0 \ \mbox{,} \\
\lim_{k \rightarrow +\infty} G(n,k) & = &
\lim_{k \rightarrow +\infty} k^{2\:b+2\:c-2\:n-4-a} \: T(n) \\
                                    & = & 0  \ \mbox{.}
\end{eqnarray*}
Again, $ \displaystyle \sum_k F(n,k) $ is independent of $n$.\\
Since $ \displaystyle \lim_{n \rightarrow +\infty} F(n,k)
= \delta_{k,0} $, we conclude again that:
\begin{equation}
\sum_{k=0}^{\infty} F(n,k)
 = \lim_{n \rightarrow +\infty} \sum_{k=0}^{\infty} F(n,k)
 = \sum_{k=0}^{\infty} \lim_{n \rightarrow +\infty} F(n,k)
 = \sum_{k=0}^{\infty} \delta_{k,0} = 1 \ \mbox{.}
\end{equation}

\subsection{A list of theorems}

These formulas can automatically be proved using the \texttt{HYPERG} package:
\begin{itemize}
\item Gauss's theorem \cite[(III.3)]{b-slater}:
      \[ {}_2F_1\left[{a,b \atop c};1\right]
         = \frac{ \Ga(c) \: \Ga(c-a-b) }
                { \Ga(c-a) \: \Ga(-b+c) } \ \mbox{.} \]

\item another Dixon's theorem \cite[(III.10)]{b-slater}:
      \[
      {}_4F_3\left[{a,1+\frac{a}{2},b,c \atop \frac{a}{2},1+a-b,1+a-c};
                                                                      -1\right]
      =
      \frac{\Ga(1+a-b) \: \Ga(1+a-c)}
           {\Ga(1+a) \: \Ga(1+a-b-c)} \ \mbox{.}
      \]

\item another Dixon's theorem \cite[(III.12)]{b-slater}:
      \begin{eqnarray*}
      & & {}_5F_4\left[{a,\frac{a}{2}+1,b,c,d \atop \frac{a}{2},1+a-b,1+a-c,1+a-d};
                                                                       1\right]
      \\
      & = &
      \frac{\Ga(1+a-b) \: \Ga(1+a-c) \: \Ga(1+a-d) \: \Ga(1+a-b-c-d)}
           {\Ga(a+1) \: \Ga(1+a-b-c) \: \Ga(1+a-b-d) \: \Ga(1+a-c-d)} \ \mbox{.}
      \end{eqnarray*}

\end{itemize}

\end{document}